\documentclass{cimart_novo}

\usepackage{adjustbox}
\usepackage{enumerate}
\usepackage[utf8]{inputenc}
\usepackage{microtype}
\usepackage{tabu}

\newtheorem{notation}{Notation}

\VOLUME{33}
\ISSUE{3}
\YEAR{2025}
\NUMBER{4}
\DOI{https://doi.org/10.46298/cm.13307}

\title{% Please, capitalize only the first word
	On $3$-generated axial algebras of Jordan type $\frac{1}{2}$
}

\author{% Please, use "Firstname Lastname" format, without abreviations
	Ravil Bildanov, Ilya Gorshkov
}

\authorinfo[% Format: please, use "F. Name"
R. Bildanov]{% Format: please, use "Department (only if 
	Sobolev Institute of Mathematics and Novosibirsk State University, Russia}{%
	ravilbildanov@gmail.com 
}

\authorinfo[% Format: please, use "S. Name"
I. Gorshkov]{% Format: please, use "Department (only if 
	Novosibirsk State Technical University, Russia}{%
	ilygor8@gmail.com
}

\abstract{%
	Axial algebras of Jordan type $\eta$ are a special type of commutative non-associative algebras. They are generated by idempotents whose adjoint operators have the minimal polynomial dividing $(x-1)x(x-\eta)$, where $\eta$ is a fixed value that is not equal to $0$ or $1$. These algebras have restrictive multiplication rules that generalize the Peirce decomposition for idempotents in Jordan algebras. 
	A universal $3$-generated algebra of Jordan type $\frac{1}{2}$  as an algebra with $4$ parameters was constructed by I. Gorshkov and A. Staroletov. Depending on the value of the parameter, the universal algebra may contain a non-trivial form radical. In this paper, we describe all semisimple $3$-generated algebras of Jordan type $\frac{1}{2}$ over a quadratically closed field.
}



\begin{document}
	
	% Here is where the main text should be typed:
	\section*{Introduction}
	Axial algebras of Jordan type $\eta$ were introduced by Hall, Rehren, and Shpectorov~\cite{hrs}
	within the framework of the general theory of axial algebras. These algebras are commutative non-associative algebras over a field $\mathbb{F}$, generated by special idempotents known as primitive axes. While Jordan algebras generated by primitive idempotents are an example of algebras of Jordan type $\frac{1}{2}$, not all algebras of this type are Jordan algebras. The Matsuo algebras, constructed from the group of $3$-transpositions, are examples of such algebras. 
	It was proved in \cite{hrs} (with a correction in \cite{hss}) that for $\eta\neq \frac{1}{2}$, algebras of Jordan type $\eta$ are the Matsuo algebras or their quotient algebras.
	Therefore, the case $\eta=\frac{1}{2}$ is special for algebras of Jordan type, and for this $\eta$, they
	are called algebras of Jordan type $\frac{1}{2}$. The class of Matsuo algebras was introduced
	by Matsuo \cite{m} and later generalized in \cite{hrs}.
	
	Algebras of Jordan type $\frac{1}{2}$ are not exhausted by Matsuo algebras and their quotient algebras. Moreover, the quotient algebras of Matsuo algebras do not contain all Jordan algebras generated by primitive idempotents. For example, the $27$-dimensional Albert algebra is generated by $4$ primitive idempotents and hence it is an algebra of Jordan type $\frac{1}{2}$ but not a Matsuo algebra \cite{kms}. 
	
	A universal $3$-generated algebra of Jordan type $\frac{1}{2}$ $A(\alpha, \beta,\gamma,\psi)$ was constructed in \cite{gs}. It is proved there that if $(\alpha+\beta+\gamma-2\psi-1)(\alpha\beta\gamma-\psi^2)\neq 0$ and $\psi^2-\alpha\beta\gamma$ is a square in $\mathbb{F}$, then $A(\alpha, \beta,\gamma,\psi)$ is isomorphic to the algebra $M^+_3(\mathbb{F})$ of $3\times 3$ matrices with Jordan multiplication. Otherwise, the algebra $A(\alpha, \beta,\gamma,\psi)$ is not simple. 
	
	A Frobenius form $(\cdot,\cdot )$ on algebra $A$ is a nonzero symmetric bilinear form that associates with multiplication in $A$, i.e., $\forall a,b,c \in A$, we have $(ab,c)=(ac,b)$ and $(a,a)\neq 0$ for each axis $a\in A$. Hall, Rehren, and Shpectorov \cite{hrs} showed that for algebras of Jordan type, there exists a unique Frobenius form with the property $(a,a)=1$ for every primitive axis $a$.
	
	Let $A$ be an algebra with a Frobenius form $(\cdot,\cdot )$. The radical of the form $(\cdot,\cdot )$ is the ideal $R(A)$ generated by elements $x$ such that $(x,a)=(a,x)=0$ for every element $a\in A$.
	
	The purpose of this article is to describe all $3$-generated algebras of Jordan type $\frac{1}{2}$ with trivial radical over a quadratically closed field. 
	
	The universal $3$-generated algebra $A(\alpha, \beta,\gamma,\psi)$ from \cite{gs} is a Jordan algebra, so any quotient algebra is also Jordan. However, we use notation from universal algebra instead of notation from Jordan algebras. 
	When we embed a 3-generated algebra into some large algebra, its parameters remain unchanged, even though the large algebra may not be Jordan.
	For example, 2-generated subalgebras are studied in \cite[Section 3]{gss}. In this paper, it is shown that isomorphic 2-generated subalgebras with different parameters of the Frobenius form have fundamentally different properties when embedded in a larger algebra. This statement emphasizes the importance of studying semisimple 3-generated algebras of Jordan type as algebras with parameters.
	
	\section{Preliminary results}
	
	We consider commutative non-associative algebras over a ground field $\mathbb{F}$ of characteristic not two. For definitions, we almost always follow~\cite{hrs} and~\cite{gg}.
	
	We denote the linear span of the set $X$ over $\mathbb{F}$ by $L \langle X \rangle$, and the algebra generated by the set $X$ by $\langle\langle X \rangle\rangle$.
	
	\begin{notation} 
		Given $a \in A$ and $\lambda \in \mathbb{F}$, consider the subspace $A_{\lambda} (a)=\{u \in A \ | \ au=\lambda u\}$.
	\end{notation}
	
	Obviously, $A_{\lambda} (a)$ is the eigenspace of the operator $ad_a: x \rightarrow ax$, associated with $\lambda \in \mathbb{F}$. 
	
	\begin{definition}
		An idempotent $a \in A$ is said to be primitive if $\dim(A_1(a)) = 1$.
	\end{definition}
	
	\begin{definition} \enlargethispage{-\baselineskip}
		An algebra $A$ is an algebra of Jordan type $\frac{1}{2}$ if $A$ is generated by a set of primitive idempotents $X$ with the following property: For every $x \in X$, there is a decomposition $A=A_0(x) \oplus A_1(x) \oplus A_{\frac{1}{2}}(x)$ with the following fusion (multiplication) rules: 
		$$A_0(x)A_{\frac{1}{2}}(x) \subseteq A_{\frac{1}{2}}(x), \quad A_1(a)A_{\frac{1}{2}}(x) \subseteq A_{\frac{1}{2}}(x),\quad  A_0(x)A_1(x) = 0,$$ $$ (A_0(x))^2 \subseteq A_0(x), \quad (A_1(x))^2 \subseteq A_1(x), \quad (A_{\frac{1}{2}}(x))^2 \subseteq A_0(x) \oplus A_1(x).$$  
		
		%\begin{center}
		%\renewcommand{\arraystretch}%{1.4}
		%\begin{tabular}{|c||c|c|c|}
		%\hline
		%$\ast$&$0$&$1$&$\frac{1}{2}$\\
		%\hline\hline
		%$0$&$0$&$ $&$\frac{1}{2}$\\
		%\hline
		%$1$&$ $&$1$&$\frac{1}{2}$\\
		%\hline
		%$\frac{1}{2}$&$\frac{1}{2}$&$\frac{1}{2}$&$\{0,1\}$\\
		%\hline
		%\end{tabular} \hspace{5pt}
		%\end{center}
		%$$A_0(a)A_{\frac{1}{2}}(a) \subseteq A_{\frac{1}{2}}(a), A_1(a)A_{\frac{1}{2}}(a) \subseteq A_{\frac{1}{2}}(a), A_0(a)A_1(a) \subseteq \{0\}$$ 
		%$$ A_0^2(a) \subseteq A_0(a), A_1^2(a) \subseteq A_1(a), A_{\frac{1}{2}}^2(a) \subseteq A_0 \oplus A_1.$$  
	\end{definition} 
	
	Such idempotents are called axes. By an $n$-generated algebra we mean an algebra generated by $n$ primitive axes. 
	
	Given an axial algebra of Jordan type $\frac{1}{2}$ $A$ and an axis $a\in A$, the map $\tau_a: A \rightarrow A$
	which acts as $x^{\tau_a}=(-1)^{2\lambda}x$, for $x\in A_{\lambda}(a)$, is an involution automorphism of $A$ called Miyamoto involution. It is clear that $y^{\tau_a}=y-4ya+4\phi_a(y)$ for each $y\in A$, where $\phi_a(y)$ is the projection of $y$ on $A_1(a)$. 
	
	Let us introduce some classes of simple Jordan algebras.  
	
	\begin{definition}
		Denote by $M_n^{+}(\mathbb{F})$ the matrix algebra $M_n(\mathbb{F})$ with Jordan product $$A \circ B = \frac{1}{2}(AB+BA).$$
	\end{definition}
	
	\begin{definition}
		If $j$ is an involutory
		automorphism of $M_n(\mathbb{F})$, then we define the Hermitian Jordan algebra $H(M_n(\mathbb{F}), j)$ as $\{ A \in M_n^{+}(\mathbb{F}) \ | \ j(A)=A \}$.  
	\end{definition}
	
	\begin{definition}
		Let $V$ be an $n$-dimensional vector space over $\mathbb{F}$ and $\phi$ a symmetric bilinear form. We define   $$ (a \oplus \mathbf{v}) \bullet (b \oplus \mathbf{w}) = (ab + \phi(\mathbf{v}, \mathbf{w})) \oplus (a\mathbf{w} + b\mathbf{v}), \text{ where } a,b \in \mathbb{F} \text{ and } \mathbf{v}, \mathbf{w} \in V.$$
		
		Then, the vector space $\mathbb{F}\oplus V$ with multiplication $\bullet$ gives a Jordan algebra which we call the Jordan algebra of the symmetric bilinear form and denote by $JForm_n(\mathbb{F})$.   
		
	\end{definition}
	It is well known that $M_n^{+}(\mathbb{F})$, $H_n^{+}(\mathbb{F})$ and $JForm_n(\mathbb{F})$ are simple Jordan algebras for $n \geq 2$ (see \cite[Theorem 23.1.2]{Mc}), so they are 
	algebras of Jordan type $\frac{1}{2}$.

	\begin{lemma}{\cite[Theorem 4.1]{hss}}\label{ff}
		Every algebra of Jordan type $\eta$ admits a unique Frobenius form which satisfies the property $(a,a)=1$ for all axes $a \in X$.
	\end{lemma}
	
	\begin{lemma}{~\cite[Proposition 2.7]{hrs}}
		The radical of Frobenius form $R(A)$ coincides with the largest ideal of $A$ containing no axes from $A$.
	\end{lemma}
	
	\begin{definition}
		An algebra of Jordan type $\eta$ is semisimple if the radical of Frobenius form is trivial.
	\end{definition}
	
	If $A$ is a Jordan type $\eta$ algebra, then $A/R(A)$ is a Jordan type $\eta$ algebra too. 
	It follows from Lemma \ref{ff} that $A/R(A)$ has a unique Frobenius form with the property $(\bar{a},\bar{a})=1$ for every axis $\bar{a}\in A/R(A)$. In the next lemma we prove that when factorizing by the radical, the Frobenius form does not change.
	
	\begin{lemma}
		Let $A$ be an algebra of Jordan type $\eta$. Then, for all $a,b \in A$ and their images $ \overline{a}, \overline{b} \in A \slash R(A)$, we have $(a,b)=( \overline{a},  \overline{b})$.
	\end{lemma}
	\begin{proof}
		Let $a=\overline{a}+r_a, b=\overline{b}+r_b$, where $ \overline{a},  \overline{b} \in A \slash R(A)$ and $ r_a, r_b \in R(A)$. 
		Then, \[(a,b)=(\overline{a}+r_a, \overline{b}+r_b)=(\overline{a}, \overline{b})+(\overline{a}, r_b)+(\overline{b}, r_a)+(r_a, r_b)=(\overline{a}, \overline{b}).
		\qedhere\]
	\end{proof}
	
	\begin{lemma}{\cite[Lemma 2]{gg}}
		Let $A$ be a finitely generated algebra of Jordan type $\frac{1}{2}$, with $a,b$ axes, and let $\alpha:=(a,b)$. Then we have the following equalities:
		\begin{enumerate}
			\item $(a_0(b))^2=(1-\alpha)a_0(b)$;
			\item $(a_{\frac{1}{2}}(b))^2=\alpha a_0(b)+(\alpha-\alpha^2)a$;
			\item $a_0(b)a_{\frac{1}{2}}(b)=\frac{1}{2}(1-\alpha)a_{\frac{1}{2}}(b)$.
		\end{enumerate}
	\end{lemma}
	
	From Lemma \ref{ff} it follows that on any algebra of Jordan type $\eta$ $A$ there is a unique Frobenius form with the property $(a,a)=1$ for any axis $a\in A$.
	
	\begin{lemma}\label{twogen} 
		Let $A=\langle\langle a,b\rangle\rangle$ be a $2$-generated algebra of Jordan type $\frac{1}{2}$.
		Then one of the following holds:
		\begin{enumerate}
			\item $\dim(A) = 1$, $(a,b)=1$, $a=b$, $\dim(R(A))=0$; 
			\item $\dim(A) = 2$, $(a,b)=0$, $A \cong \mathbb{F} \oplus \mathbb{F}$, $\dim(R(A))=0$;
			\item $\dim(A) = 2$, $(a,b)=1$, $\dim(R(A))=1$;
			\item $\dim(A) = 3$, $(a,b)=0$, $\dim(R(A))=1$, $A/R(A) \cong \mathbb{F} \oplus \mathbb{F}$;
			\item $\dim(A) = 3$, $(a,b)=1$, $\dim(R(A))=2$;
			\item $\dim(A) = 3$, $(a,b) \neq 0,1$, and $A$ is a Matsuo algebra. In particular, it is a simple Jordan algebra isomorphic to $JForm_2(\mathbb{F})$.  
		\end{enumerate}
	\end{lemma}
	\begin{proof}
		The assertion of the lemma is a simple consequence of \cite[Theorem 1.1]{hrs}. 
	\end{proof}
	
	\begin{lemma}{\cite[Corollary 1]{gg}} 
		Let $A$ be a $2$-generated algebra of Jordan type $\frac{1}{2}$ with generating axes $a$ and $b$. Let $\alpha:=(a,b)$. Then we have
		\begin{enumerate}
			\item $a(ab)=\frac{1}{2}(\alpha a+ab)$;
			\item $(ab)b=\frac{1}{2}(\alpha b+ab)$;
			\item $(ab)(ab)=\frac{\alpha}{4}(a+b+2ab)$.
		\end{enumerate}
	\end{lemma}
	
	\begin{lemma}{\cite[Theorem 1]{gs}}\label{mainGS}
		Let $A$ be a $3$-generated Jordan type $\frac{1}{2}$ algebra.
		There exists a $3$-generated $9$-dimensional algebra $A(\alpha,\beta,\gamma,\psi)$ such that $A$ is a quotient algebra of $A(\alpha,\beta,\gamma,\psi)$ for suitable values of parameters $\alpha,\beta,\gamma,\psi$. 
	\end{lemma}
	
	Let $A=\langle\langle a,b,c \rangle\rangle$, $\dim(A)=9$, $\alpha=(a,b)$, $\beta=(b,c)$, $\gamma=(a,c)$, $\psi=(ab,c)$. In Table~1 below (that is similar to~\cite[Table 6]{gs} up to renumbering rows), we present all possible relations for $\alpha,\beta,\gamma,\psi$ for $A(\alpha,\beta,\gamma,\psi)$ to not be simple. 
	
	\begin{table}[h!]\label{tb}
		\begin{center}
			%			\thisfloatpagestyle{empty}
			%\begin{adjustbox}{angle=90}
			\begingroup
			\setlength{\tabcolsep}{20pt} %
			\renewcommand{\arraystretch}{1.3}
			
			{\tiny
				$\begin{tabu}[h!]{|c|c|c|c|c|c|c|}
					\hline
					\text{$A_i$}&\text{Relations} & \dim(A/R(A)) & \text{Basis of the radical}  \\ \hline
					
					A_1&\psi=\alpha=\beta=\gamma=1 & 1 & 
					\begin{tabu}{@{}c@{}}
						b-a, c-a, ab-a, bc-a, ac-a, \\ a(bc)-a, b(ac)-a, c(ab)-a
						
					\end{tabu}
					\\ \hline
					%case 4
					A_2&\psi=\alpha=\beta=0, \gamma=1 & 2 & 
					\begin{tabu}{@{}c@{}} 
						c-a, ab, bc, ac-a, \\ a(bc), b(ac), c(ab)
						
					\end{tabu}
					\\ \hline
					
					%case 3
					A_3&\psi=\alpha=\beta=\gamma=0 & 3 & 
					\begin{tabu}{@{}c@{}} 
						ab, bc, ac, a(bc), b(ac), c(ab)
					\end{tabu}
					\\ \hline
					%case 5
					A_4&\begin{tabu}{@{}c@{}}
						\psi=\alpha=0, \beta, \gamma\neq0, \\
						\beta+\gamma=1 
					\end{tabu}
					& 3 & 
					\begin{tabu}{@{}c@{}}
						
						ab, \cfrac{1}{2}\gamma a-\cfrac{1}{2}\beta b-\cfrac{1}{2}c+bc, \\ -\cfrac{1}{2}\gamma a+\cfrac{1}{2}\beta b-\cfrac{1}{2}c+ac, \\ \cfrac{1}{4}\gamma a+\cfrac{1}{4}\beta b-\cfrac{1}{4}c+a(bc), \\
						\cfrac{1}{4}\gamma a+\cfrac{1}{4}\beta b-\cfrac{1}{4}c+b(ac), c(ab)   
					\end{tabu}
					\\ \hline
					
					%case 2 
					A_5&\begin{tabu}{@{}c@{}} \alpha\beta\gamma=\psi^2, \psi\neq0,\alpha\neq1, \\ \alpha+\beta+\gamma=2\psi+1 \end{tabu} & 3 &
					\begin{tabu}{@{}c@{}} 
						
						\alpha(\beta-1)a+\alpha(\gamma-1)b+\alpha(1-\alpha)c+(2\alpha-2\psi)ab, \\ (\alpha\beta-\alpha\psi)b+(\psi-\alpha\beta)ab+(\alpha^2-\alpha)bc, \\ (\alpha\gamma-\alpha\psi)a+(\psi-\alpha\gamma)ab+(\alpha^2-\alpha)ac, \\ (\alpha\psi-\alpha^2\beta)a+(\alpha+\psi-\alpha^2-\alpha\gamma)ab+2\alpha(\alpha-1)a(bc), \\ \alpha(\psi-\alpha\gamma)b+(\alpha+\psi-\alpha^2-\alpha\beta)ab+2\alpha(\alpha-1)b(ac), \\ (\psi-\alpha\beta)a+(\psi-\alpha\gamma)b+(1-\alpha)ab+2(\alpha-1)c(ab)
					\end{tabu}  \\ \hline
					
					%case 7
					A_6&\psi=\alpha=\beta=0, \gamma\neq0,1 & 4 &
					\begin{tabu}{@{}c@{}} 
						ab, bc, ac, a(bc), b(ac), c(ab)
					\end{tabu} 
					\\ \hline
					
					%case 8
					A_7&\begin{tabu}{@{}c@{}}
						\psi^2\neq\alpha\beta\gamma, \\
						\alpha+\beta+\gamma=2\psi+1, \\ \alpha\neq1 
					\end{tabu} & 4 &
					\begin{tabu}{@{}c@{}}
						\frac{1}{2}(\beta-1)a+\frac{1}{2}(\beta-\alpha)b+\frac{1}{2}(1-\alpha)c+(1-\beta)ab+(\alpha-1)bc, \\ \frac{1}{2}(\gamma-\alpha)a+\frac{1}{2}(\gamma-1)b+\frac{1}{2}(1-\alpha)c+(1-\gamma)ab+(\alpha-1)ac, \\ (2\psi-2\alpha\beta+\beta-1)a+(\gamma-1)b+(1-\alpha)c+(4-2\alpha-2\gamma)ab+(4\alpha-4)a(bc), \\ (\beta-1)a+(2\psi-2\alpha\gamma+\gamma-1)b+(1-\alpha)c+(4-2\alpha-2\beta)ab+(4\alpha-4)b(ac), \\ (\psi-\alpha)a+(\psi-\alpha)b+\alpha(1-\alpha)c+(2-\beta-\gamma)ab+(2\alpha-2)c(ab) 
					\end{tabu}
					\\ \hline
					%case 6
					A_8&\begin{tabu}{@{}c@{}}
						\psi=\alpha=0, \beta,\gamma\neq0, \\
						\beta+\gamma\neq1 
					\end{tabu}
					& 6 & 
					\begin{tabu}{@{}c@{}}
						
						ab, b(ac)-a(bc), c(ab)
					\end{tabu} \\ \hline
					%case 1
					A_9&\begin{tabu}{@{}c@{}} \alpha\beta\gamma=\psi^2, \psi\neq0, \\ \alpha+\beta+\gamma\neq2\psi+1 \end{tabu} &  6 &
					\begin{tabu}{@{}c@{}}
						-\beta\gamma ab-\alpha\beta ac+2\psi a(bc),\\ -\beta\gamma ab-\alpha\gamma bc+2\psi b(ac), \\ -\alpha\gamma bc-\alpha\beta ac+2\psi c(ab)
						
					\end{tabu}    \\ \hline
				\end{tabu}$ } 
			
			\caption{Bases of the radical}\label{t:radical}
			\endgroup
			%\end{adjustbox}
		\end{center}
		%\end{sidewaystable}
	\end{table}
	
	\section{Main Results}
	
	In this section we assume that $A$ is a $3$-generated algebra of Jordan type $\frac{1}{2}$ with a trivial radical over a quadratically closed field $\mathbb{F}$ and denote by $(\cdot,\cdot)$ the unique Frobenius form on $A$ satisfying the property that $(a,a)=1$ for every axis $a$ of $A$.
	
	\begin{theorem}
		Let $A$ be a $3$-generated algebra of Jordan type $\frac{1}{2}$ with a trivial radical over a quadratically closed field $\mathbb{F}$ with characteristic not equal to two or three. Then A is isomorphic to one of the following algebras:
		\enlargethispage{-\baselineskip}
		\begin{enumerate}
			\item $\mathbb{F}^n, n \in \{ 1, 2, 3 \}$;
			\item $JForm_2(\mathbb{F})$;
			\item $\mathbb{F} \oplus JForm_2(\mathbb{F})$;
			\item $M^{+}_2(\mathbb{F})$;
			\item $H(M_3(\mathbb{F}), j)$ with $j(X)=X^T$;
			\item $M^{+}_3(\mathbb{F})$.
		\end{enumerate}
	\end{theorem}
	
	It follows from Lemma \ref{mainGS}, that any algebra of Jordan type $\frac{1}{2}$ is isomorphic to the quotient algebra of $A(\alpha, \beta, \gamma, \psi)$ for some parameters $\alpha, \beta, \gamma, \psi$. From \cite[Enlightenment Structure Theorem]{Mc} we obtain that $A(\alpha, \beta, \gamma, \psi)/R(A(\alpha, \beta, \gamma, \psi))$ is a direct sum of simple algebras. Consequently, to describe 3-generated algebras of Jordan type $\frac{1}{2}$ with a trivial radical, we need to describe the quotient algebras of $A(\alpha, \beta, \gamma, \psi)$ by its radical for each choice of the parameters $\alpha, \beta,\gamma, \psi$.
	
	%It follows from Lemma \ref{mainGS}, that we need to describe the quotient algebras of the algebra $A(\alpha, \beta, \gamma, \psi)$ by its radical.
	We use the description of the algebra $A(\alpha, \beta, \gamma, \psi)$ from \cite[Theorem 2]{gs}.
	Let us recall that, following \cite{gs}, we use the notation $\alpha=(a,b), \beta=(b,c), \gamma=(a,c), \psi=(ab,c)$. 
	
	In Table \ref{tb}, one can find the dimensions and bases of the radicals of the algebra $A(\alpha, \beta, \gamma, \psi)$. Denote by $A_i$ the universal $9$-dimensional algebra $A(\alpha_i, \beta_i, \gamma_i, \psi_i)$ with parameters and numeration from Table $1$, $R_i$ the radical of this algebra and by $S_i$ the quotient algebra $A_i/R_i$. 
	
	We begin with two trivial propositions for $1$-dimensional and $2$-dimensional algebras, which are not generated by three linearly independent axes.  
	
	\begin{proposition}
		If $A$ is a $1$-dimensional algebra of Jordan type $\frac{1}{2}$ with a trivial radical, then $A \cong S_1$. 
	\end{proposition}
	\begin{proof}
		It is easy to see that $S_1 \cong \mathbb{F}$. 
		We have that $A$ is 1-dimensional, so $\dim L\langle a,b,c\rangle=1$ and $a=b=c$. Hence $A \cong \mathbb{F} \cong S_1$.    
	\end{proof}

	\begin{proposition}
		If $A$ is a $2$-dimensional $3$-generated algebra of Jordan type $\frac{1}{2}$ with a trivial radical, then $A \cong \mathbb{F} \oplus \mathbb{F}\cong S_2$. 
	\end{proposition}
	\begin{proof}
		By Lemma \ref{twogen}, there is only one $2$-dimensional algebra of Jordan type $\frac{1}{2}$ with a trivial radical, so $A \cong \mathbb{F} \oplus \mathbb{F} \cong S_2$. 
	\end{proof}

	\begin{lemma}
		Algebras $S_4$ and $S_5$ are isomorphic.
	\end{lemma} 
	\begin{proof}
		
		We will first show that $S_4=\langle\langle a,c \rangle\rangle$. Put $S=\langle\langle a,c \rangle\rangle$. We have that $$S=S_{0}(a)+S_{1}(a)+S_{\frac{1}{2}}(a) \text{ and } c=c_0(a)+\gamma a+c_{\frac{1}{2}}(a),$$ where $c_0(a) \in S_0(a)$ and $c_{\frac{1}{2}}(a) \in S_{\frac{1}{2}}(a)$. 
		Firstly, assume that $c_{\frac{1}{2}}(a)=0$. In this case, we have that $c=c^2=(c_0(a))^2+(\gamma a)^2$, in particular $(c_0(a))^2=c_0(a)$. Consequently, $c_0(a)c=\nolinebreak(c_0(a))^2=\nolinebreak c_0(a)$; which contradicts $c$ being the primitive idempotent. Thus, $c_{\frac{1}{2}}(a) \neq 0$. Assume that $c_0(a)=0$. Then, $$\gamma a+c_{\frac{1}{2}}(a)=c=c^2=\gamma^2 a+\gamma c_{\frac{1}{2}}(a)+(c_{\frac{1}{2}}(a))^2.$$ 
		Hence, $\gamma=1$ and from the definition of $S_4$ it follows that $\beta=0$. In this case, we have that $(a-c,b)=(a,b)-(c,b)=0$. It follows that $a-c \in R(S_4)$, which is a contradiction. Therefore, $\dim (S)=3$. Thus, $S_4=S$ and $S_4$ is generated by $2$ axes. From Lemma \ref{twogen} it follows that $S_4\simeq JForm_2(\mathbb{F})$
		
		Now, consider the algebra $S_5$. By definition of $S_5$, we have $(a,b)\not\in\{0,1\}$, where $a,b$ are the axes from the generating set of the algebra $S_5$. Lemma \ref{twogen} implies that $\langle\langle a,b\rangle\rangle\simeq JForm_2(\mathbb{F})$. Since $dim(\langle\langle a,b\rangle\rangle)=3$, we obtain $\langle\langle a,b\rangle\rangle=S_5$, in particular $S_5\simeq S_4$.
	\end{proof}

	\begin{proposition}
		If $A$ is a $3$-dimensional $3$-generated algebra of Jordan type $\frac{1}{2}$ with a trivial radical, then $A$ is isomorphic to either $S_3$ or $S_5$.  
	\end{proposition}
	
	\begin{proof}
		Assume that $A$ is 2-generated and let $a$ and $b$ are the generating axes. It follows from Lemma~\ref{twogen} that there is only one $3$-dimensional $2$-generated algebra of Jordan type $\frac{1}{2}$ with a trivial radical. In this case, we can choose any other axis of the algebra $A$. Let us call it $c$. Put $c=a^{\tau_b}=a-4ab+4\alpha b$. We have $\beta=\alpha,\ \gamma=(1-2\alpha)^2$ and $\psi=\alpha(2\alpha-1)$. Therefore $\alpha\beta\gamma=\psi^2,\ \psi\neq0,\alpha\neq1,$ and $ \alpha+\beta+\gamma=2\psi+1$. So, in this case, $A\simeq S_5$.  
		
		Assume that $A$ is not generated by $2$ axes. Therefore, based on the dimension of $A$, we conclude that $A=L\langle a, b, c\rangle$.
		
		Assume that $ab \notin L\langle a,b \rangle$. Hence $\dim \langle\langle a, b \rangle\rangle = 3$. Therefore $c \in L\langle a,b,ab \rangle=A$, which is a contradiction.
		Similarly, we can show that $ac \in L\langle a,c \rangle$ and $bc \in L \langle b, c \rangle$. In particular, we have $\dim(\langle\langle a,b\rangle\rangle)=\dim(\langle\langle a,c\rangle\rangle)=\dim(\langle\langle c,b\rangle\rangle)=2$. From Lemma \ref{twogen} it follows that $\{(a,b),(a,c),(b,c)\}\subseteq \{0,1\}$. Moreover, if $(a,b)=0$, then $\langle\langle a,b\rangle\rangle\simeq \mathbb{F}\oplus\mathbb{F}$. Therefore, if $(a,b)=(a,c)=(b,c)=0$, then $A\simeq \mathbb{F}\oplus\mathbb{F}\oplus\mathbb{F}$ and $\psi=0$. In this case the Gram matrix of the algebra $A$ is the identity matrix and hence the radical of $A$ is trivial.
		We conclude that in this case $A\simeq S_3$.
		
		Assume that $(a,c)\neq0$. We have $(a,c)=1$. In this case, $R(\langle\langle a,c\rangle\rangle)$ is not trivial and contains the element $a-c$. Assume that $(a,b)=(b,c)=0$. In this case we have $(a-c,b)=0$. Consequently, $a-c\in R(A)$, which is a contradiction. Therefore, without loss of generality, we can assume that $(b,c)=1$. If $(a,b)=1$ then $(a-c,b)=0$ and consequently $a-c\in R(A)$, which is a contradiction. Therefore $(a,b)=0$. From the description of $2$-generated algebras of Jordan type $\frac{1}{2}$ we have $ab=0$, $a= c+a_h$, $b=c+ b_h$, where $a_h,b_h\in A_{\frac{1}{2}}(c)$. Therefore, $$0=ab=(c+a_h)(c+b_h)=c+\frac{1}{2}(a_h+b_h)+a_hb_h,$$ where $c+a_hb_h\in A_{0}(c)\oplus A_{1}(c)$ and $a_h+b_h\in A_{\frac{1}{2}}(c)$. Therefore $a_h+b_h=0$. In particular, $b=a^{\tau_c}$ and $\dim(A)=2$.  
		
		Lemma \ref{twogen} implies that in this case $A\simeq JForm_2(\mathbb{F})$.
	\end{proof}

	\begin{lemma}
		$S_6$ is isomorphic to $\mathbb{F} \oplus JForm_2(\mathbb{F})$.
	\end{lemma}
	\begin{proof}
		Let $\langle\langle a,b,c\rangle\rangle\simeq S_6$. We have $(a,c)\not\in\{0,1\}$. Therefore $\langle\langle a,c\rangle\rangle$ is isomorphic to $JForm_2(\mathbb{F})$. From Table 1, it follows that the radical of $A(0,0,\gamma,0)$ contains $ab$ and $bc$. Therefore, $ab=bc=0$ and $S_6\simeq \langle\langle a,c\rangle\rangle\oplus \langle\langle b\rangle\rangle\simeq \mathbb{F} \oplus JForm_2(\mathbb{F})$.  
	\end{proof}
	\begin{lemma}
		$S_7$ is isomorphic to $M^{+}_2(\mathbb{F})$. 
	\end{lemma}
	\begin{proof} 
		
		Assume that $\alpha=\beta=\gamma=0$. We have $(c^{\tau_b},a)=2\alpha\beta+\gamma-4\psi=-2$ and $\langle\langle c^{\tau_b},a,b\rangle\rangle=\langle\langle
		a,b,c\rangle\rangle=S_7$. Thus, up to redesignation of the generating elements, we can assume that $\gamma\neq0$.
		
		Let 
		\begin{equation*}
			A = \left(
			\begin{array}{ccc}
				1-\lambda_c & 1 \\ \lambda_c(1-\lambda_c) & \lambda_c \\
			\end{array}
			\right),
			B = \left(
			\begin{array}{ccc}
				1 & 0 \\
				\lambda_b & 0 \\
			\end{array}
			\right), 
			C =\left(
			\begin{array}{cc}
				1 & \lambda_a \\
				0 & 0 \\
			\end{array}
			\right),  
		\end{equation*}
		where $\lambda_a,\lambda_b, \lambda_c\in \mathbb{F}\setminus \{0\}$. Consider the following map $f: S_7 \rightarrow M^{+}_2(\mathbb{F}),$ defined by $ f(a)=A, f(b)=B, f(c)=C$.
		It is easy to see that $\dim L\langle A, B, C, A \circ B \rangle=4$, so $\langle\langle A, B, C \rangle\rangle=M^{+}_2(\mathbb{F})$.
		
		A map $( \cdot \ , \cdot ): M^{+}_2(\mathbb{F})^2 \rightarrow \mathbb{F}$ such that $(X,Y)=tr(XY)=tr(X \circ Y)$, for $X,Y\in M^{+}_2(\mathbb{F})$, is a symmetric bilinear form on $M^{+}_2(\mathbb{F})$ (see \cite[Chapter 1.6]{Mc}). This form associates with the product $\circ$. Clearly, we have $tr(A \circ A)=tr(B \circ B)=tr(C \circ C)=1$. 
		
		Furthermore, we see that %$tr(A \circ B)=1 - \lambda_a + \lambda_b=\alpha, tr(B \circ C)=1+\lambda_a\lambda_b=\beta, tr(A \circ C)=1-\lambda_c+\lambda_a\lambda_c(1-\lambda_c)=\gamma$ and $tr(A \circ (B \circ C))=tr(B \circ (A \circ C))=tr(C \circ (A \circ B))=\psi=\frac{1}{2}(1-\alpha-\beta-\gamma)$. 
		\begin{align*}
			tr(A \circ B)&=1 - \lambda_a + \lambda_b=\alpha, \\
			tr(B \circ C)&=1+\lambda_a\lambda_b=\beta, \\
			tr(A \circ C)&=1-\lambda_c+\lambda_a\lambda_c(1-\lambda_c)=\gamma \text{ and} \\
			tr(A \circ (B \circ C))&=tr(B \circ (A \circ C))=tr(C \circ (A \circ B))=\psi=\frac{1}{2}(1-\alpha-\beta-\gamma).
		\end{align*}
		If $\alpha\neq0$ then: 
		\begin{align*}
			\lambda_a&=\frac{1}{\alpha(\alpha-1)}(\psi+\alpha\gamma \mp \sqrt{\psi^2-\alpha\beta\gamma}), \\
			\lambda_b&=-\frac{-1 + \beta + \gamma}{\gamma}, \\
			\lambda_c&=\frac{1-\beta}{\gamma}
		\end{align*}
		%$$\lambda_a=\frac{1}{\alpha(\alpha-1)}(\psi+\alpha\gamma \mp \sqrt{\psi^2-\alpha\beta\gamma}), \lambda_b=\frac{1}{\gamma}(\psi+\alpha\gamma \pm \sqrt{\psi^2-\alpha\beta\gamma}),$$ $$\lambda_c=\pm \frac{1}{\gamma}(\psi+\alpha+\sqrt{\psi^2-\alpha\beta\gamma}),$$ 
		are the solution of these equations.
		
		If  $\alpha=0$ then: $$ \lambda_a=-\frac{\gamma^2-\gamma}{-1 + \beta +\gamma},\quad \lambda_b=-\frac{-1 + \beta + \gamma}{\gamma}, \quad \lambda_c=\frac{1-\beta}{\gamma}$$ are the solution of these equations. Note that in this case $-1+\beta+\gamma\neq0$, otherwise $\psi=\alpha\beta\gamma=0$.
		
		Using computer calculations, we show that multiplication table for $f(\langle\langle a,b,c \rangle\rangle)$ \footnote{Computer calculations for multiplication table in $S_7$ can be found in \url{https://github.com/RavilBildanov/3gen-axial-algebras/blob/main/S7 multiplication table.nb}, see paragraph Tables.} coincides with multiplication table for $S_7$. Hence, $f$ is an isomorphism. 
		
		We also use computer calculations to check that $R(f(\langle\langle a,b,c \rangle\rangle))=\{0\}$ and relations between $\alpha, \beta, \gamma, \psi$ hold \footnote{One can find our computer calculations here: \\ \url{https://github.com/RavilBildanov/3gen-axial-algebras/blob/main/M2+ (S7).nb}}.
	\end{proof}
	
	% \begin{table}[h]
		% \begingroup
		% \begin{center}
			% \footnotesize
			% \begin{tabu}{| c || c | c | c | c |}
				% \hline
				% $\ast$  & $a$ & $b$  & $c$ & $ab$ \\ \hline\hline
				% $a$  & $a$  & *  & * & * \\ \hline
				% $b$  & $ab$ & $b$ & * & *  \\ \hline
				% $c$  & \begin{tabu}{@{}c@{}}$\frac{1}{2(\alpha-1)}((\gamma-\alpha)a+$\\
					% $+(\gamma-1)b+(1-\alpha)c+$\\
					% $+2(-\gamma+1)ab)$\end{tabu} & \begin{tabu}{@{}c@{}}$\frac{1}{2(\alpha-1)}((\beta-1)a+(\beta-\alpha)b$ \\$+(1-\alpha)c+2(-\beta+1)ab)$\end{tabu}  & $c$ & * \\ \hline
				% $ab$ &  $\frac{1}{2}(a\alpha+ab)$ & $\frac{1}{2}(b\alpha+ab)$ & 
				%         \begin{tabu}{@{}c@{}}$\frac{1}{2(\alpha-1)}((\psi-\alpha)a+$ \\$+(\psi-\alpha)b+(\alpha-\alpha^2)c+$\\$+(2-\beta-\gamma)ab)$\end{tabu} & 
				%         \begin{tabu}{@{}c@{}}$\frac{1}{4}\alpha(a+$ \\ $+b+2ab)$\end{tabu} \\ \hline
				% \end{tabu}
			% \caption{Multiplication table for $S_7$}\label{t:prod_S7}
			% \end{center}
		% \endgroup
		% \end{table}
	
	\begin{table}[h]
		\renewcommand{\arraystretch}{1.5}
		\centering
		\footnotesize
		\setlength{\tabcolsep}{3pt}
		\begin{tabu}{| c || c | c | c | c |}
			\hline
			$\ast$  & $a$ & $b$  & $c$ & $ab$ \\ \hline\hline
			$a$  & $a$  & *  & * & * \\ \hline
			$b$  & $ab$ & $b$ & * & *  \\ \hline
			$c$  & $\begin{array}{l}\frac{1}{2(\alpha-1)}\big((\gamma-\alpha)a\\
				\:\:+(\gamma-1)b+(1-\alpha)c\\
				\:\:+2(-\gamma+1)ab\big)\end{array}$ & $\begin{array}{l}\frac{1}{2(\alpha-1)}\big((\beta-1)a+(\beta-\alpha)b \\\:\:+(1-\alpha)c+2(-\beta+1)ab\big)\end{array}$  & $c$ & * \\ \hline
			$ab$ &  $\frac{1}{2}(a\alpha+ab)$ & $\frac{1}{2}(b\alpha+ab)$ & 
			$\begin{array}{l}\frac{1}{2(\alpha-1)}\big((\psi-\alpha)a \\\quad +(\psi-\alpha)b+(\alpha-\alpha^2)c\\\quad+(2-\beta-\gamma)ab\big)\end{array}$ & 
			$\begin{array}{l}\frac{1}{4}\alpha(a+b \\ \quad\quad+2ab)\end{array}$ \\ \hline
		\end{tabu}
		\caption{Multiplication table for $S_7$}\label{t:prod_S7}
	\end{table}

	\begin{proposition}
		If $A$ is a $4$-dimensional $3$-generated algebra of Jordan type $\frac{1}{2}$ with a trivial radical, then one of the following assertions holds: 
		
		\begin{enumerate}
			\item $A \simeq S_6 \simeq \mathbb{F} \oplus JForm_2(\mathbb{F})$;
			\item $A \simeq S_7 \simeq M^{+}_2(\mathbb{F})$.
		\end{enumerate}
	\end{proposition}
	
	\begin{proof}
		The algebra $M^{+}_2(\mathbb{F})$ is a simple Jordan algebra (see \cite[Chapter 1.11]{Mc}). The algebra $\mathbb{F} \oplus JForm_2(\mathbb{F})$ contains non-trivial ideals. Therefore $M^{+}_2(\mathbb{F})\not\simeq \mathbb{F} \oplus JForm_2(\mathbb{F})$. Hence, to prove this proposition, it suffices to show that $S_6 \simeq \mathbb{F} \oplus JForm_2(\mathbb{F})$ and $S_7 \simeq M^{+}_2(\mathbb{F})$. Thus the Proposition follows from Lemmas 8 and 9. 
	\end{proof}

	\begin{lemma}
		$S_8 \cong H(M_3(\mathbb{F}), j)$.
	\end{lemma}
	\begin{proof}
		Consider the following matrices in $H(M_3(\mathbb{F}), j)$ and the map $f: S_8 \rightarrow H(M_3(\mathbb{F}), j)$ defined by  $f(a)=A$, $f(b)=B$, $f(c)=C$, where
		\begin{equation*}
			A = \left(
			\begin{array}{ccc}
				1 & 0 & 0 \\
				0 & 0 & 0 \\
				0 & 0 & 0 \\
			\end{array}
			\right)
			B = \left(
			\begin{array}{ccc}
				0 & 0 & 0 \\
				0 & \frac{1+\sqrt{1-4\lambda_b^2}}{2} & \lambda_b \\
				0 & \lambda_b & \frac{1-\sqrt{1-4\lambda_b^2}}{2} \\
			\end{array}
			\right)
			C = \left(
			\begin{array}{ccc}
				\frac{1+\sqrt{1-4\lambda_c^2}}{2} & 0 & \lambda_c \\
				0 & 0 & 0 \\
				\lambda_c & 0 & \frac{1-\sqrt{1-4\lambda_c^2}}{2} \\
			\end{array}
			\right).
		\end{equation*}
		The scalars $\lambda_a, \lambda_b, \lambda_c \in \mathbb{F}\setminus\{0\}$ are the invariant by $\theta$ parameters which are defined later from conditions imposed on $\alpha,\beta,\gamma$ and $\psi$.

		We show that the mapping $f$ is an isomorphism between the algebras $S_8$ and $H(M_3(\mathbb{F}), j)$. It is easy to see that $A^2=A$, $B^2=B$, and $C^2=C$. We check that 
		\[f(\langle\langle a,b,c \rangle\rangle)=L\langle A, B, C, A \circ C, B \circ C, A \circ (B \circ C) \rangle.\]
		Thus, $\dim L\langle A, B, C, A \circ C, B \circ C, A \circ (B \circ C) \rangle=6$. Hence $\langle\langle A,B,C \rangle\rangle$ and $H(M_3(\mathbb{F}), j)$ are isomorphic as vector spaces.   
		
		A map $( \cdot, \cdot ): H(M_3(\mathbb{F}), j)^2 \rightarrow \mathbb{F}$ such that $(X,Y)=tr(XY)=tr(X \circ Y)$, where $X,Y\in H(M_3(\mathbb{F}), j)$ is a symmetric bilinear form on $H(M_3(\mathbb{F}), j)$ \cite[Chapter 1.6]{Mc}. This form associates with the product $\circ$. Clearly, we have $tr(A \circ A)=tr(B \circ B)=tr(C \circ C)=1$. 
		Furthermore, we see that
		\begin{align*}
			tr(A \circ B)&=0, \\
			tr(B \circ C)&=\frac{1}{4}(1-\sqrt{1-4\lambda_b^2})(1-\sqrt{1-4\lambda_c^2})=\beta, \\
			tr(A \circ C)&=\frac{1+\sqrt{1-4\lambda_c^2}}{2}=\gamma, \text{ and } \\
			tr(A \circ (B \circ C))&=tr(B \circ (A \circ C))=tr(C \circ (A \circ B))=0.
		\end{align*}
		% $tr(A \circ B)=0, tr(B \circ C)=\frac{1}{4}(1-\sqrt{1-4\lambda_b^2})(1-\sqrt{1-4\lambda_c^2})=\beta, tr(A \circ C)=\frac{1+\sqrt{1-4\lambda_c^2}}{2}=\gamma$ and $tr(A \circ (B \circ C))=tr(B \circ (A \circ C))=tr(C \circ (A \circ B))=0$.
		So, we have conditions on $\lambda_b, \lambda_c$.
		
		Take the basis $\{a, \:b, \:c, \:b \cdot c, \:a \cdot c, \:a \cdot (b \cdot c)\}$ for $S_8$. The multiplication table for $f(\langle\langle a,b,c \rangle\rangle)$ coincides with multiplication table for $S_8$\footnote{Computer calculations for the multiplication table in $S_8$ can be found in \url{https://github.com/RavilBildanov/3gen-axial-algebras/blob/main/S8 multiplication table.nb}, see paragraph Tables.}. % We can build multiplication table for $S_8$ by the same way as the multiplication table for $S_7$. 
		
		We also use computer calculations to check that $R(f(\langle\langle a,b,c \rangle\rangle))=\{0\}$ and relations between $\alpha, \beta, \gamma, \psi$ hold \footnote{Computer calculations for this proof can be found in \\ \url{https://github.com/RavilBildanov/3gen-axial-algebras/blob/main/H3+ (S8).nb}}.
		
		We have that $\{A, B, C, B \circ C, A \circ C, A \circ (B \circ C)\}$ is a basis of the algebra $H(M_3(\mathbb{F}), j)$ and hence the kernel of $f$ is trivial. Thus $f$ is an isomorphism of the algebras $S_8$ and $H(M_3(\mathbb{F}), j)$.
	\end{proof}

	\begin{lemma}
		Algebras $S_8$ and $S_9$ are isomorphic.
	\end{lemma}
	\begin{proof}
		
		Assume that $\alpha=\beta=\gamma=1$, but then $\psi=1$ and we obtain a contradiction. Let $\alpha \neq 1$, take $d=x_a(b)=\dfrac{2ab-\alpha a-b}{\alpha-1}$. It is known from \cite{gg}, that $d$ is a primitive idempotent in $S_9$ with $ad=0$ and so $d$ is an axis because $S_9$ is a Jordan algebra. 
		
		Assume that $(c,d)\neq0$. It can be proved via computer calculations that the set $\{a,c,d,ac,cd,a(cd)\}$ is an additive basis of $B=\langle\langle a,c,d\rangle\rangle$. In particular, $B=S_9$. Define a homomorphism $f$ from $S_9$ to $S_{8}$ given by $f(a)=\bar{a}, f(d)=\bar{b}, f(c)=\bar{c}$. We have $(\bar{a},\bar{b})=0$, $(\bar{a}, \bar{b}\bar{c})=0$, $(\bar{a}, \bar{c})=\bar{\beta} \neq 0$ and $(\bar{b},\bar{c})=\bar{\gamma}\neq 0$. The relation $\bar{\beta}+\bar{\gamma}=\frac{2\psi-\beta-\alpha\gamma}{\alpha-1}+\gamma \neq 1$ is equivalent to $\alpha+\beta+\gamma \neq 2\psi+1$, so $f$ is an isomorphism between $S_8$ and $S_9$ iff $(c,d) \neq 0$.    
		
		Now assume that $(c,d)=0$. We will prove that $\gamma \neq 1$ in this case. We have that 
		\begin{align*}
			0&=(c, (\alpha-1)d)=(c,2ab-\alpha a-b, c) \\ &
			=(2ab-b,c)-\alpha=(2ab,c)-(b,c)-\alpha \\
			&=2\psi-\beta-\alpha.
		\end{align*}
		%$$0=(c, (\alpha-1)d)=(c,2ab-\alpha a-b, c)=$$ $$=(2ab-b,c)-\alpha=(2ab,c)-(b,c)-\alpha=2\psi-\beta-\alpha$$
		If $\gamma=1$, then $2\psi+1=\beta+\alpha+\gamma$, a contradiction. Therefore, $\gamma \neq 1$. 
		
		Put $d'=\dfrac{2ac-\gamma a-c}{\gamma-1}$. Note that $(b,d') \neq 0$. Indeed, if $(b,d')=0$, we have 
		\begin{align*}
			0&=(b, 2ac-\gamma a-c)=(b, 2ac)-(b, \gamma a)-(b,c)\\ &=2\psi-\gamma\alpha-\beta
		\end{align*}
		%$$ 0=(b, 2ac-\gamma a-c)=(b, 2ac)-(b, \gamma a)-(b,c)=2\psi-\gamma\alpha-\beta$$
		
		Using these equalities, we obtain $\alpha=\gamma\alpha$ and so $\gamma=1$ or $\alpha=0$, a contradiction.
		
		We then take $a,b,d'$ as the new generating set of $S_9$ and, using the same computer calculations, prove that $\{a,b,x,ab,ax,b(ax)\}$ is an additive basis of $B'=\langle\langle a,b,d' \rangle\rangle $. Define a homomorphism $f'$ from $S_9$ to $S_{8}$ given by $f(a)=\bar{a}'$, $f(d')=\bar{b}'$, and $f(b)=\bar{c}'$. We have $(\bar{a}',\bar{b}')=0$, $(\bar{a}',\bar{b}'\bar{c}')=0$, $(\bar{a}', \bar{c}')=\bar{\beta}' \neq 0$ and $(\bar{b}',\bar{c}')=\bar{\gamma}'\neq 0$. The relation $\bar{\beta}'+\bar{\gamma}'=\frac{2\psi-\beta-\alpha\gamma}{\gamma-1}+\alpha \neq 1$ is equivalent to $\alpha+\beta+\gamma \neq 2\psi+1$, so $f'$ is an isomorphism between $S_8$ and $S_9$ iff $(b,d') \neq 0$.
		
		We can also check that the multiplication table for $a,c,d,ac,cd,a(cd)$ coincides with the multiplication table for the standard basis $a,b,c,bc,ac,a(bc)$ of $S_8$. This means that $S_9$ contains a 6-dimensional subalgebra isomorphic to $S_8$ \footnote{One can see our computer calculations here:  \url{https://github.com/RavilBildanov/3gen-axial-algebras/blob/main/3gen axial algebra.nb}, see section "Isomorphism between $S_8$ and $S_9$".}.  
	\end{proof}
	
	\begin{table}[h]
		\renewcommand{\arraystretch}{1.5}
		\centering
		\footnotesize
		\setlength{\tabcolsep}{3pt}
		\begin{tabu}{|c||c|c|c|c|c|c|}
			\hline
			$\ast$   & $a$ & $b$ & $c$ & $bc$ & $ac$ & $a(bc)$ \\ \hline\hline
			$a$  & $a$  & * & * & *  & *  & *  \\ \hline
			$b$  &  0 & $b$  & * & *  & *  & *  \\ \hline
			$c$  &  $ac$  & $bc$  & $c$  & *  & *  & *  \\ \hline
			$bc$ &  $a(bc)$ & $\frac{1}{2}(b\beta+bc)$  &  $\frac{1}{2}(c\beta+bc)$  &  $\begin{array}{l}\frac{\beta}{4}(b+c+2bc)\end{array}$ & *  & *  \\ \hline
			$ac$ &  $\frac{1}{2}(\gamma a+ac)$ & $a(bc)$  &  $\frac{1}{2}(\gamma c+ac)$ &  $\begin{array}{l}\frac{\gamma}{4}bc + \frac{\beta}{4}ac \\
				\quad+\frac{1}{2}a(bc)\end{array}$  &  $\begin{array}{l}\frac{\gamma}{4}(a+c\\
				\quad+2ac)\end{array}$  & *  \\ \hline
			$a(bc)$ & $0$  & $\begin{array}{l}\frac{1}{4}(\beta ac+2a(bc))\end{array}$  & $\begin{array}{l}\frac{1}{4}(\gamma bc + \beta ac)\end{array}$  & $\begin{array}{l}\frac{\beta\gamma}{8}b+\frac{\beta}{8}ac\\
				\quad +\frac{\beta}{4}a(bc)\end{array}$  & $\begin{array}{l}\frac{\beta\gamma}{8}a+\frac{\gamma}{8}bc\\
				\quad+\frac{\gamma}{4}a(bc)\end{array}$  &  $\begin{array}{l}\frac{\beta\gamma}{16}(a 
				+b)\end{array}$  \\ \hline
		\end{tabu}
		\caption{Multiplication table for $S_8$}\label{t:prod_S9}
	\end{table}
	
	\begin{proposition}\label{S8}
		If $A$ is a $6$-dimensional $3$-generated algebra of Jordan type $\frac{1}{2}$ with a trivial radical, then $A \simeq S_8 \simeq S_9\simeq H(M_3(\mathbb{F}), j)$, where $j(X)=X^T$.
	\end{proposition}
	\begin{proof}
		The proposition follows from Lemmas 10 and 11.
	\end{proof}
	
	\subsection*{Acknowledgements}
	The work was supported by the Theoretical Physics and Mathematics
	Advancement Foundation "BASIS". We are grateful to V. A. Afanasev, who read the draft of this paper and pointed out several errors.

	%%% REFERENCES %%%
	%{\small\bibliography{Ilya}}
	% Please, do not change the above line and do not insert your references
	% into this file.  Instead, insert your references into the cimart.bib file.
	% See cimart.bib for further instructions.

	\EditInfo{March 28, 2024}{June 7, 2024}{Ivan Kaygorodov}
	
\end{document}